\documentclass{article}
\usepackage[latin1]{inputenc}
\usepackage{amsmath}
\usepackage{amssymb}
\usepackage{amsthm}
\usepackage{mathrsfs}

\begin{document}

\theoremstyle{plain} \newtheorem{Thm}{Theorem}
\newtheorem{prop}{Proposition}
\newtheorem{lem}{Lemma}
\newtheorem{de}{Definition}
\newtheorem{rem}{Remark}
\theoremstyle{remark} \newtheorem*{pf}{Proof}
\renewcommand\theenumi{(\alph{enumi})}
\renewcommand\labelenumi{\theenumi}
\renewcommand{\qedsymbol}{}
\renewcommand{\qedsymbol}{\ensuremath{\blacksquare}}
\newcommand{\stb}{, \ldots ,}
\newcommand{\stp}{+ \ldots +}
\title{Non-expander Cayley graphs of simple groups}
\author{G\'abor Somlai \\ Alfr\'ed R\'enyi Institute of Mathematics \\ email: zsomlei@gmail.com\thanks{Research is partially supported by MTA Renyi "Lendulet" Groups and Graphs Research Group}}
\date{}
\maketitle

\begin{abstract}
For every infinite sequence of simple groups of Lie type of growing rank we exhibit connected Cayley graphs of degree at most $10$ such that the isoperimetric number of these graphs converges to $0$. This proves that these graphs do not form a family of expanders.
\end{abstract}

\section{Introduction}
Let $G$ be a finite group and $T$ a subset of $G$. The Cayley graph
$Cay(G,T)$ is defined by having vertex set $G$ and $g$ is adjacent to
$h$ if and only if $g^{-1} h \in T$. A Cayley graph $Cay(G,T)$ is undirected if
and only if $T= T^{-1}$, where $T^{-1} = \left\{ t^{-1} \in G \mid t
  \in T \right\}$.

Let $\Gamma$ be an arbitrary graph and $S \subseteq V(\Gamma)$. We
define the boundary of $S$ which we denote by $\partial S$ to be the
set of vertices in $V(\Gamma) \setminus S$ with at least one neighbour
in $S$. For a graph $\Gamma$ the isoperimetric number $h(\Gamma)$ is defined by
\[ h(\Gamma) = min \left\{ \frac{\left| \partial S
  \right| }{\left| S \right| } \mbox{ } \bigg| \mbox{ } S \subset
V(\Gamma) \mbox{, } 0 < \left| S \right| \le \frac{ \left| V(\Gamma)
  \right| }{2} \right\} \mbox{.} \]
A graph $\Gamma$ is called an $\epsilon$-expander if $h(\Gamma) \ge \epsilon$ and a
series of $k$-regular graphs $\Gamma _n$ is called an expander family if there is a
constant $\epsilon >0$ such that for every $n$ the graph $\Gamma _n$ is an $\epsilon$-expander.
Finally, we say that a family of groups $G_n$ is a family of uniformly expanding groups if there exist $0 < \epsilon \in \mathbb{R}$
and $k \in \mathbb{N}$ such that for every $i$ and every generating set $S_i \subset G_i$ of size at most $k$ the Cayley graphs $Cay(G_i, S_i)$ are $\epsilon$-expanders.

The study of series of Cayley graphs of finite simple groups has received great attention.
The proof of the fact was completed by Kassabov, Lubotzky and Nikolov in \cite{3} based on several earlier work (see \cite{10},\cite{11},\cite{12},\cite{13},\cite{14}), that there exist $k \in \mathbb{N}$ and $0 < \epsilon \in \mathbb{R}$ such that every
non-abelian finite simple group which is not a Suzuki group has a set of generators $S$ of size at most $k$ for which $Cay(G,S)$ is an $\epsilon$-expander.
This work was extended by Breuillard, Green and Tao in \cite{4} to the Suzuki groups. These results can also motivate the question which was asked by Lubotzky in \cite{5} whether every family of Chevalley groups of bounded rank is a family of uniformly expanding group.

Lubotzky also suggested to investigate families of simple groups of unbounded rank in \cite{2} and wrote that it seems likely that if
$G_n$ is a sequence of non-abelian simple groups such that the rank of $G_n$ is unbounded, then for every $n$ there exists a generating set $T_n \subset G_n$ such that the graphs $Cay(G_n, T_n)$ do not form a family of expanders. An explicit example (see \cite{7}) of a non-expander family of Cayley graphs of special linear groups was given by Luz. The diameter of the graphs given by Luz was investigated by Kassabov an Riley and it was proved in \cite{8} that there exists $c \in \mathbb{R}$ such that the diameter of the graphs is smaller than $c \mbox{ }log (\left| SL(n,p)\right|)$.
Similarly, the symmetric group $S_n$ is generated by $\gamma =(12)$ and $\sigma_n = (1,2, \ldots , n)$ for every $n \in \mathbb{N}$ and the sequence of isoperimetric numbers $h\left( Cay(S_n, \gamma, \sigma_n) \right)$ tends to $0$, see \cite{7}. Moreover, one can find a set of generators of $S_n$ such that the diameter of the corresponding Cayley graphs is $\Omega (n^2)$ which gives that these Cayley graphs do not form a family of expanders.

We will investigate $7$ series ($A_l, B_l, C_l, D_l, A_{2n-1}^{1}, A_{2n}
^{1}, D_{n}^{1} )$ of finite simple
groups of Lie type. These are the groups of Lie type such that the rank of a sequence of groups tends to infinity if we fix a series of the Lie group.
In order to define generators and subgroups of these groups we will
use the generators given by Steinberg in \cite{1} and we will use the notation and several results of the book of Carter \cite{9}.

For these $7$ series of finite
simple groups of Lie type we construct Cayley graphs and subsets
such that the number of the neighbours of these subsets depends on the
rank of the groups. Moreover, the isoperimetric number of these graphs tends
to $0$. This proves the conjecture of Lubotzky concerning the series of Cayley graphs of simple groups of unbounded rank.
More precisely, we prove the following:
\begin{Thm}\label{thm1}
\begin{enumerate}
\item\label{ta} Let $G$ be a Chevalley group of rank $l$ of type $A_l$, $B_l$, $C_l$ or $D_l$.
For every $l \ge 5$ and for every finite field $GF(q)$ there exists a generating set $T$ of cardinality at most $10$ and a subset of the vertices $S \subset V(Cay(G,T))$ with $\left| S \right| \le \frac{\left| G \right|}{2}$ such that $\frac{\left| \partial(S) \right|}{\left| S
  \right|} \le \frac{6}{ l-3}$.
\item\label{tb} Let $G$ be a twisted group of type $A_{2n-1}^{1}$, $D_{n}^{1}$ or $A_{2n}
^{1}$. For every $n \ge 5$ and for every finite field $GF(q)$ there exists a generating set $T'$ of cardinality at most $8$ and $S' \subset V(Cay(G,T'))$ with $\left| S' \right| \le \frac{\left| G \right|}{2}$ such that $\frac{\left| \partial(S') \right|}{\left| S'
\right|} \le \frac{6}{n-2}$.
\end{enumerate}
\end{Thm}

The paper is organized into the following 4 sections. In Section \ref{pre} we give all necessary definitions and we collect some important fact about the construction of simple groups of Lie type. The proof of Theorem \ref{thm1} \ref{ta} is contained in Section \ref{chevalley} and Theorem \ref{thm1} \ref{tb} which is the case of twisted groups will be handled in Section \ref{twisted}. In Section \ref{identification} we present the original construction in terms of matrices which was extended to several different series of simple groups.

\section{Preliminaries}\label{pre}
In this section we collect important facts about about finite simple
groups of Lie-type and we build up the notation we will use all along this paper.

Let $K=GF(q)$ be a finite field.
We denote by $\Phi$ the system of roots and $\Phi = \Phi ^{+} \cup
\Phi^{-}$ is the union of the positive and negative roots. We also
choose $\Pi =\{r_1, r_2, \ldots , r_l \}\subset \Phi^{+}$ which is the set of the fundamental roots.

The Weyl group $W$ is generated by the fundamental reflections $w_r$,
where $r \in \Pi$. In order to simplify notation we denote by $w_i$ the fundamental
reflections $w_{r_i}$, where $r_i \in \Pi$. We denote by $x_r(\psi)$ the standard generators of the Chevalley group $G$, where
$r \in \Phi$ and $\psi \in K$. If $r= r_i$ for some $r_i \in \Pi$, then we denote by $x_i(\psi)$ the standard generator $x_{r_i}(\psi)$.
These elements generate the Chevalley
group $G$. The subgroups $X_{r} = \left\{ x_r (t) \mid  t \in K \right\}$ are called
root subgroups of $G$ if $r \in \Phi$.

The Weyl group $W$ is isomorphic to
$N/H$ for some $H \lhd N \le G$. The cosets of $H$ in $N$ can be written
as $n_w H$ for all $w \in W$ and $N$ is generated by $H$ and the
elements $n_r$ for $r \in \Phi$. Moreover, $n_r =x_r(1) x_{-r}(-1)x_r(1)$ is the element of the subgroup generated by the root subgroups $X_r$ and $X_{-r}$.
It is well known that $n_r x_s (t) n_r ^{-1}
=x_{w_r(s)}(\eta_{r,s} t )$ for some $\eta_{r,s} \in K$ depending only on $r$ and $s$, see \cite[p.101.]{9}.
The elements of the normal subgroup $H$ of $N$ can be written in the form $h(\chi)$
where $\chi$ is a $K$-character of $\mathbb{Z} \Phi$. The subgroup $H$ is generated by the elements
of the set $\left\{ h_r (\lambda) \mid r \in \Phi \mbox{ and } \lambda \in K^{*} \right\}$, where the $K$-character corresponding
to $h_r (\lambda)$ is $\chi_{r, \lambda}$ with $\chi_{r, \lambda} (a) = \lambda ^{\frac{2(a,r)}{(r,r)} }$.
$H$ is a normal subgroup of $N$ and $n_w h (\chi) n_w^{-1} = h (\chi ')$, where $\chi '(r) = \chi (w^{-1} (r))$ see \cite[p.102.]{9}.
Furthermore, $h_r(\lambda) = n_r(\lambda) n_r(-1)$ and hence $h_r(\lambda) \in \langle X_r, X_{-r} \rangle$, see \cite[p.96.]{9}.
\section{Chevalley groups}\label{chevalley}
In this section we construct series of Cayley graphs for $4$ different
series of Chevalley groups. For these Chevalley groups we need $6$ series of Cayley graphs. The six different constructions are similar
but we will treat them separately.

We first prove the following technical lemma.
\begin{lem}\label{lem}
Let $w=w_1 w_2 \ldots w_l$ be a Coxeter element of the Weyl group $W$
and let us assume that the fundamental root $r_i$ is orthogonal to $r_j$ if $i+1 < j \le l$ and $r_{i+1}$ is orthogonal to $r_k$ if $1 \le k \le i-1$. We also assume that $r_i$ and $r_{i+1}$ have the same length and $w_i (r_{i+1}) =r_i +r_{i+1}$. Then $w(r_i) = r_{i+1}$.
\end{lem}
\begin{pf}
Since $r_i$ is orthogonal to $r_j$ for every $j>i+1$ we have that $w(r_i)= w_1 w_2 \ldots w_i w_{i+1} (r_i)$. The elements $w_k$ are reflections through the hyperplane perpendicular to $r_k$. Thus $w_k (r_k)= -r_k$ for every $1 \le k \le l$ and $w_{i+1} (r_i) =r_i +r_{i+1} =w_{i} (r_{i+1} )$ since $r_i$ and $r_{i+1}$ have the same length.
It follows that $w_i w_{i+1} (r_i)= w_i (r_i +r_{i+1}) =w_i(r_i) + w_{r_{i+1} } =-r_i + (r_i +r_{i+1}) =r_{i+1}$. Hence $w(r_i) =w_1 w_2 \ldots w_{i-1} (r_{i+1}) = r_{i+1}$ since $r_{i+1}$ is orthogonal to $r_k$ for $1 \le k \le i-1$.
\qed
\end{pf}

\subsection{$A_l$}\label{A_l}
Let $G$ be a Chevalley group of type $A_l$.
The Dynkin diagram of the corresponding root system is the following.

\unitlength 1mm
\begin{picture}(100,10)
\put(5,5){ \circle{1} }
\put(5,7){$r_1$ }
\put(7,5){ \line(1,0){8} }

\put(19,5){ \circle{1} }
\put(19,7){$r_2$ }
\put(21,5){ \line(1,0){8} }

\put(33,5){ \circle{1} }
\put(33,7){$r_3$ }
\put(35,5){ \line(1,0){8} }

\put(45,5){ \mbox{\dots }}
\put(60,5){ \circle{1} }
\put(60,7){$r_{l-1}$ }
\put(63,5){ \line(1,0){8} }
\put(58,5){ \line(-1,0){8} }
\put(74,5){ \circle{1} }
\put(74,7){$r_l$ }

\end{picture}

One can see from the Dynkin diagram that  $w_i (r_{i+1} )= r_i +r_{i+1} = w_{i+1} (r_{i} )$ for $i=1, \ldots ,l-1$.

Let $w=w_1 w_2 \ldots w_l$ be a Coxeter element of the
Weyl group. We choose $\lambda$ to be a generator of the
multiplicative group of $GF(q)$.
\begin{lem}\label{gena}
$x_1(1)$, $n_w$ and $h_{r_1} (\lambda )$ generate the Chevalley group $G$.
\end{lem}
\begin{pf}

It was proved in \cite{1} that $x_1 (1) n_w$ and $h_{r_1} (\lambda)$ generate $G$.
Clearly, $x_1(1)$ and $n_w$ generate $x_1(1) n_w$ which proves the Lemma.
\qed
\end{pf}

For every $l\ge 5$ we define the following undirected Cayley graph:
\[ \Gamma_a = Cay(G,\left\{ x_1(1), n_w, h_{r_1}(\lambda),  x_1(1)^{-1},
  n_w^{-1}, h_{r_1}(\lambda) ^{-1} \right\} ) \mbox{.} \]

Let $K_a$ be the subgroup of the Chevalley group $G$ generated by the root subgroups $X_{r_1}, X_{-r_{1} },  X_2,
X_{-r_{2} },  \ldots , X_{r_{l-1}}, X_{-r_{l-1 }}$
and let \[ S_a = \cup_{i=0}^{l-1} K_a n_w^{i} \mbox{.} \]
Every element of the Weyl group $W$ acts on the the root system
$\Phi$.
\begin{lem}\label{orbita}
The orbit of $w$ which contains $r_1$ is the following:

\unitlength 1mm
\begin{picture}(110,10)

\put(0,5){\vector(1,0){4} }

\put(5,4,5){$r_1$ }
\put(7,5){ \vector(1,0){8} }
\put(11,6){\mbox{w}}
\put(19,4,5){$r_2$ }
\put(21,5){ \vector(1,0){8} }
\put(25,6){\mbox{w}}

\put(33,4,5){$r_3$ }
\put(35,5){ \vector(1,0){6} }
\put(38,6){\mbox{w}}
\put(43,5){ \mbox{\dots }}

\put(48.5,5){ \vector(1,0){6} }
\put(51.5,6){\mbox{w}}

\put(57,4,5){$r_{l-1}$ }
\put(62,5){ \vector(1,0){8} }
\put(66,6){\mbox{w}}

\put(73,4.5){$r_l$ }
\put(76,5){ \vector(1,0){8} }
\put(80,6){\mbox{w}}

\put(86,5){$-r_1 - \ldots - r_{l}$ }
\put(106,5){ \vector(1,0){8} }
\put(110,6){\mbox{w}}
\end{picture}

This can be formulated as follows: \begin{equation*} \begin{split}
w(r_i) &=r_{i+1} \mbox{ }for \mbox{ }1 \le i \le l-1 \\ w(r_l) &= -r_1 - \ldots - r_l \\
 w(-r_1 - r_2 -\ldots - r_l ) &= r_1 \end{split} \end{equation*}
\end{lem}
\begin{pf}
Lemma \ref{lem} gives that $w(r_i) = r_{i+1}$  for $1 \le i \le n-1$ and
\[ w(r_l) =w_1 w_2 \ldots w_l (r_l ) =w_1 w_2 \ldots w_{l-1} (-r_l ) = -w_1 w_2 \ldots w_{l-1} (r_l )\]
since $w$ is a linear transformation of the vector space spanned by the roots.
We also have $w_j(r_{j+1} + \dots + r_l) = r_j + r_{j+1} + \dots + r_l$ for $ 1 \le j \le l-1$. Therefore
\begin{equation*} \begin{split} w_1 w_2 \ldots w_{l-1} (r_l ) &= w_1 w_2 \ldots w_{l-2} (r_{l-1}+r_l )\\ &= w_1 w_2 \ldots w_{l-3} (r_{l-2}+r_{l-1} +r_l ) = \dots = r_1 + r_2 + \ldots + r_l \mbox{.} \end{split} \end{equation*}
This shows that \begin{equation}\label{eq2} w (r_l) = - (r_1 +r_2 + \ldots + r_l) \mbox{.} \end{equation}
\\Using again the linearity of $w$ and equation (\ref{eq2}) we get \[ w(r_1 +r_2 + \ldots +r_l ) =r_2 +r_3 \stp r_l - \left( r_1 + r_2 + \ldots + r_l \right) = -r_1 \mbox{.} \]
This gives $w (-(r_1 +r_2 \stp r_l ) ) =r_1$, finishing the proof Lemma \ref{orbita}.
\qed
\end{pf}
It follows from Lemma \ref{orbita} that if $1 \le i \le l-1$, then $n_w^{i} K_a n_w^{-i}$
contains $n_w^{i} X_{r_{l-i}} n_w^{-i} =X_{r_l}$. Therefore $n_w^{i} K_a n_w^{-i}
\ne K_a$ which shows that $n_w^{i} \notin K_a$ for every $1 \le i \le
l-1$. This implies that $K_a, K_a n_w \stb K_a n_w^{l-1}$ are different
right cosets of $K_a$ so $S_a$ is the union of $l$ pairwise disjoint subsets of the vertices of $\Gamma_a$
and these subsets have the same cardinality.
\begin{lem}
  $ \frac{\left| \partial(S_a) \right|}{\left| S_a \right| } \le \frac{6}{l}$
\end{lem}
\begin{pf}
$S_a$ is the union of $l$ right cosets of $K_a$ so $\left| S_a \right|
= l \left| K_a \right| $.
It is clear from the definition of $S_a$ that $(K_a n_w^{i}) n_w
\subset S_a$ for every $0 \le i \le l-2$ and similarly $\left( K_a n_w^{i} \right)
n_w^{-1} \subset S_a$ if $1 \le i \le l-1$. Therefore those neighbors of $S_a$
which are not in $S_a$ can only be obtained as an element of following subset of
the vertices of $\Gamma_a$: \begin{equation*} \begin{split} &K_a n_w^{l} \bigcup K_a n_w^ {-1} \bigcup_{i=1}^{l-1} \left( K_a n_w^{i} \right)
x_1(1) \bigcup_{i=1}^{l-1} \left( K_a n_w^{i} \right) x_1(1)^{-1} \bigcup_{i=1}^{l-1} \left( K_a n_w^{i}
\right) h_r(\lambda) \\ &\bigcup_{i=1}^{l-1} \left( K_a w^{i} \right)
h_r(\lambda)^{-1} \mbox{.} \end{split} \end{equation*}

$K_a$ is a subgroup of $G$ so $\left( K_a n_w^{i} \right) x =K_a n_w^{i}$ if and only if $n_w^{i}
x n_w^{-i} \in K_a$. We first apply this observation to $x_1(1)$ and $x_1
(1)^{-1} = x_1(-1)$. It is easy to see from Lemma \ref{orbita} that $n_w^{i} x_1(\pm 1) n_w^{-i}$ is of the form $x_{w^{i} \left( r_1 \right) } (\alpha) = x_{i+1}(\alpha)$ for some
$\alpha \in GF(q)^{*}$ if $0 \le i \le l-1$.
It follows that  $n_w^{i} x_1(\pm 1) n_w^{-i} \in X_{r_{i+1}}  \subset K_a$ if $i \ne l-1$.

Using the fact that $h_r (\lambda)$ and $h_r (\lambda)^{-1}$ are in the subgroup $\langle X_r, X_{-r} \rangle$ we get that $n_w^{i} h_{r_1}
   (\lambda )^{\pm 1} n_w^{-i} \in \langle X_{w^{i}(r_1)  }, X_{-w^{i}(r_1) } \rangle= \langle X_{r_{i+1} }, X_{-r_{i+1} } \rangle \subset K_a$
  if $i \ne l-1$.

Now, $\partial S_a \subseteq K_a n_w^l \cup K_a n_w^{-1} \cup K_a n_w^{l-1} x_1(1) \cup K_a n_w^{l-1} x_1 (-1) \cup
K_a n_w^{l-1} h_{r_1} (\lambda) \cup K_a n_w^{l-1} h_{r_1} \left( \frac{1}{\lambda}
\right)$. These subsets are all of them right cosets of $K_a$ so they
have the same cardinality which proves that $\left| \partial S_a
\right| \le 6 \left| K_a \right| $, while $\left| S_a \right| = l
\left| K_a \right|$.
\qed

\end{pf}

\begin{rem}
In order to prove Theorem \ref{thm1} \ref{ta} we repeat the previous construction several times. In every single case the connection set of the Cayley graph will consist of few standard generators of the Chevalley group, an element of the form $n_w$, where $w= w_1 w_2 \ldots w_l$ is a Coxeter element of the corresponding Weyl group and an element of the group $H$.
If $G$ is of rank $l$ we will choose a subgroup of $G$ which is isomorphic to a Chevalley group of rank $l-1$ and which is of the same type. The subset of the vertices for which the isoperimetric number is sufficiently small will be the union of cosets of the subgroup of rank $l-1$.
\end{rem}

\subsection{$B_l$}
Let $G$ be a Chevalley group of type $B_l$
The Dynkin diagram of a Chevalley group of type $B_l$ is the
following:

\begin{center}
\unitlength 1mm
\begin{picture}(72,11)

\put(6,5){ \circle{2} }
\put(6,7){$r_1$ }
\put(8,5.5){ \line(1,0){8} }
\put(9,4.5){\line(1,0){8}}

\put(19,5){ \circle{2} }
\put(19,7){$r_2$ }
\put(21,5){ \line(1,0){8} }

\put(33,5){ \circle{2} }
\put(33,7){$r_3$ }
\put(35,5){ \line(1,0){6} }

\put(43,5){ \mbox{\dots }}
\put(56,5){ \circle{2} }
\put(56,7){$r_{l-1}$ }
\put(54,5){ \line(-1,0){6} }
\put(59,5){ \line(1,0){8} }
\put(69,5){ \circle{2} }
\put(69,7){$r_l$ }

\end{picture}
\end{center}
It is easy to see from the Dynkin diagram that $w_1(r_2)=r_2 + 2 r_1$ and $w_2(r_1)=r_1 + r_2$.

One can see using Lemma \ref{lem} that \begin{equation}\label{eq1} w(r_i)= w_1 w_2 \ldots w_l (r_{i}) =r_{i+1} \mbox{ } \mbox{ for } 2  \le i \le l-1. \end{equation}  The fundamental roots $r_3, \ldots , r_l$ are orthogonal to $r_1$. Therefore $w(r_1)=w_1 w_2 (r_1) =w_1 (r_1 + r_2)= -r_1 +(r_2 + 2 r_1 ) = r_1 +r_2$.
We also have that $w$ is linear so using equation (\ref{eq1}) we have that if $2 \le j \le l-1$, then
\begin{equation}\label{pot1} w(r_1 +r_2 + \ldots + r_j)=w(r_1 )+ w(r_2) + \ldots + w(r_j) = r_1 +r_2 + r_3 + \ldots + r_{j+1} \mbox{.} \end{equation}
Using these observations we conclude that the following picture represents a part of the orbit of the action of the group generated by $w$ including the root $r_1$:

\unitlength 1mm
\begin{picture}(100,10)\label{picb}

\put(0,4.5){$r_1$ }
\put(2,5){ \vector(1,0){6} }
\put(5,6){\mbox{w}}
\put(10,4.5){${r_1 + r_2}$ }
\put(20,5){ \vector(1,0){6} }
\put(23,6){\mbox{w}}
\put(28,4.5){$r_1 + r_2 +r_3$ }
\put(46,5){ \vector(1,0){6} }
\put(49,6){\mbox{w}}
\put(54,5){ \mbox{\dots }}
\put(59,5){ \vector(1,0){6} }
\put(67,4.5){$r_1 +r_2 + \ldots +r_l$ }
\put(94,5){ \vector(1,0){6} }
\put(97,6){\mbox{w}}

\end{picture}

This can be formulated as follows:
\begin{equation}\label{negy}
w^{i}(r_1)= r_1+ r_2 + \ldots + r_{i+1} \mbox{ } \mbox{for } i=1, \ldots ,l-1 \mbox{.}
\end{equation}

The orbit of $\langle w \rangle$ containing these elements contains $w(r_1 +r_2 + \ldots +
r_l)$ as well. It is easy to see that $w_i (r_{i+1} + r_{i+2} + \ldots + r_l) = r_{i} + r_{i+1} + \ldots +r_l$ if $2 \le i \le l-1$.
We also have $w_1(r_2) =2r_1+r_2$ hence
\begin{equation*} \begin{split} w(r_l) &= w_1 \ldots w_{l-1} w_l(r_l) =-w_1 \ldots w_{l-1}(r_l) \\ &=-w_1 \ldots w_{l-2}(r_{l-1}+r_l) = \dots = -w_1 (r_2  \stp r_l ) \\ & = -(r_l \stp r_2 +2r_1) \mbox{.} \end{split} \end{equation*}
This implies using equation (\ref{pot1}) that \begin{equation}\label{eq3} w(r_1 \stp r_l) =w(r_1 \stp r_{l-1} ) + w(r_l) =-r_1 \mbox{.} \end{equation}
One can easily describe the remaining elements of the orbit since $w$ is linear.

We also investigate the action of $\langle w \rangle$ on $2r_1 +r_2 + \ldots + r_l$ and $r_2 + \ldots + r_l$.
Using equation (\ref{eq3}) and the linearity of $w$ we get that $w (2 r_1 + r_2 + \ldots + r_l) = w(r_1) + w (r_1 + r_2 + \ldots + r_l)= r_1 + r_2 - r_1= r_2$.
 It follows using equation (\ref{eq1}) that \begin{equation}\label{eq4} w^{i} (2 r_1 + r_2 + \ldots + r_l) = r_{i+1} \mbox{ for } 1 \le i \le l-1 \mbox{.} \end{equation}

One can also prove using equation (\ref{negy}) and equation (\ref{eq4}) that
\begin{equation}\label{eq5}
 w^{i} (r_2 + \ldots + r_l) = -2(r_1+ r_2 \stp r_{i+1}) + r_{i+1} \mbox{ } \mbox{ for } 1 \le i \le l-1 \mbox{.}
\end{equation}
\subsubsection{$Char(K)>2$}
Let us assume that $char(K) >2$.
\begin{lem}\label{genb}
$x_1(1)$, $n_w$ and $h_t (\lambda )$, where $t=2 r_1 +r_2 + \dots + r_l$ generate the Chevalley
group $G$ of type $B_l$ if the characteristic of the underlying field is not $2$.
\end{lem}
\begin{pf}
It was proved in \cite{1} that $x_1(1) n_w$ and $h_t (\lambda )$ generate the Chevalley group $G$ if $char(K) \ne 2$.
\qed
\end{pf}

We define again a sequence of connected Cayley graphs.
Let
\[ \Gamma _b = Cay \left( G, \left\{ x_1(1), x_1(-1), n_w, n_w^{-1}, h_t
  (\lambda), h_t (\lambda) ^{-1} \right\} \right) \mbox{,} \]
  where $G$ is of rank $l$ and $w=w_1 w_2 \ldots w_l$.
 Similarly to the previous case let
\[ K_b = \langle X_{r_1}, X_{-r_{1} }, X_{r_2},
   X_{-r_{2} } \stb X_{r_{l-1} } ,  X_{-r_{l-1 }}   \rangle  \] and let
 \[ S_b = \cup_{i=0}^{l-2} K_b n_w^{i} \mbox{.} \]

\begin{lem}\label{b_1}
$ \frac{\left| \partial(S_b) \right|}{\left| S_b \right| } \le \frac{4}{l-1}$
\end{lem}
\begin{pf}
We claim that $S_b$ is the union of pairwise disjoint right cosets
of the same subgroup $K_b$ in $G$. We only have to show that $n_w^{i} \notin K_b$ if
$1 \le i \le l-2$. Straightforward calculation shows using equation (\ref{eq1}) that $n_w^{i}
X_{r_{l-i} } n_w^{-i} = X_{r_l}$ if $1 \le i \le l-2$. Therefore $ X_{r_l}
\subset n_w^{i} K_b n_w^{-i} \ne K_b$ if $1 \le i \le l-2$ which gives that $n_w^{i} \notin K_b$. Thus $S_b$ is the union of $l-1$ pairwise disjoint right cosets of $K_b$.

Using the definition of the Cayley graph $\Gamma _b$ we have that $\partial S_b$ is a subset of the following set:

\begin{equation*} \begin{split} &\bigcup\limits_{i=0}^{l-2} (K_b n_{w}^{i}) n_w \bigcup\limits_{i=0}^{l-2} (K_b n_w ^{i}) n_w ^{-1} \bigcup\limits_{i=0}^{l-2} (K_b n_w^{i}) x_1(1) \bigcup\limits_{i=0}^{l-2} (K_b n_w^{i}) x_1(-1) \\ &\bigcup\limits_{i=0}^{l-2} (K_b n_{w}^{i}) h_t ( \lambda) \bigcup\limits_{i=0}^{l-2} (K_b n_{w}^{i}) h_t ( \lambda)^{-1} \mbox{.}  \end{split} \end{equation*}
By the definition of $S_b$ the subsets $K_b n_w^{i} n_w$ are contained in
$S_b$ if $0 \le i \le l-3$ and $K_b n_w^{i} n_w^{-1} \subset S_b$ if $1 \le i \le l-2$.

Using equation (\ref{negy}) we have $n_w ^{i} x_1(\pm 1) n_w ^{-i} =
x_{r_1 + r_2 + \ldots + r_{i+1} } (t)$ for some $t \in K^{*}$. If $0 \le i \le l-2$, then $x_{r_1 + r_2 + \ldots + r_{i+1} }(t) \in K_b$ since $r_1 + r_2 + \ldots + r_{i+1}$ is in the root system generated by the fundamental roots $r_1, r_2 , \ldots , r_{l-1}$ and $K_b$ is the Chevalley group of type $B_{l-1}$ generated by the corresponding root subgroups. Therefore $K_b n_w^{i} x_1 (\pm 1) =K_b n_w^{i} \subset S_b$ if $0 \le i \le l-2$.

The elements $h_t (\lambda)$ and $h_t (\lambda)^{-1}=h_t (\frac{1}{\lambda} )$ are in the subgroup generated by $X_t$  and  $X_{-t} $.
Equation (\ref{eq4}) shows that $n_w ^{i} X_t n_w ^{-i} =X_{w^{i}(t) } =X_{r_{i+1} }$ and by the linearity of $w$ we have $n_w ^{i} X_{-r} n_w ^{-i} =X_{-r_{i+1 } }$ for $i=1,2, \ldots, l-2$. Thus $n_w ^{i} h_t (\lambda) n_w ^{-i}$ and $n_w ^{i} h_t
(\frac{1}{\lambda}) n_w ^{-i}$ are in $\langle X_{r_{i+1} }, X_{-r_{i+1} } \rangle \le K_b$ if $1 \le i \le l-2$.

It follows that $\partial S_b \subset K_b n_w^{l-1} \cup K_b n_w ^{-1 } \cup K_b h_t(\lambda) \cup K_b h_t(\frac{1}{\lambda})$ which gives $\frac{\left| \partial S_b
  \right|}{\left| S_b \right| } \le \frac{4 |K_b|}{(l-1) |K_b|} =\frac{4}{l-1}$.
  \qed
\end{pf}
\subsubsection{$Char(K) =2$}

\begin{lem}
$x_s(1)$, $x_{-r_{1}} (1)$, $n_w$ and $h_t (\lambda )$, where $s=r_2 + \dots + r_l$ generate the Chevalley
group $G$ of type $B_l$ if $char (K)= 2$.
\end{lem}
\begin{pf}
It was proved in \cite{1} that $x_s(1) x_{-r_1} (1) n_w$ and $h_t (\lambda )$ generate $G$ if $K = GF(2^k)$  with $k>1$ and
$x_s(1) x_{-r_1} (1)$ and $n_w$ generate $G$ if $\left| K \right| =2$.
\qed
\end{pf}
Let
\[ \Gamma _b' = Cay \left( G, \left\{ x_s(1), x_{-r_1}(1), n_w^{\pm 1}, h_t
  (\lambda)^{\pm 1} \right\} \right) \mbox{.} \]
The set $S_b$ can be considered as a subset of $V(\Gamma_b')$ so we claim the following.

\begin{lem}
$ \frac{\left| \partial(S_b) \right|}{\left| S_b \right| } \le \frac{5}{l-1}$
\end{lem}
\begin{pf}
It was proved in Lemma \ref{b_1} that $|S_b| = (l-1) |K_b|$.

Similarly, the proof of Lemma \ref{b_1} shows that $K_b n_w^{i} h_t(\lambda)^{\pm 1} \subset S_b$ if $1 \le i \le l-2$.
By the definition of $S_b$ we have $K_b n_w^{i} n_w \subset S_b$ if $0 \le i \le l-3$ and $K_b n_w^{i} n_w^{-1} \subset S_b$ if $1 \le i \le l-2$.

Using $w(-r_1) =-w(r_1)$ and equation (\ref{negy}) we get that $n_w^{i} x_{-r_{1} } (1) n_w^{-i} \in K_b$ if $0 \le i \le l-2$ since $w^{i} (r_1) = r_1+ r_2 + \ldots + r_{i+1} $ by equation (\ref{negy}). Hence $K_b n_w^{i} x_{-r_{1}}(1) = K_b n_w^{i} \subset S_b$.

Equation (\ref{eq5}) shows that $n_w^{i} x_s(1) n_w^{-i}$ is in $K_b$ if $1 \le i \le l-2$. Therefore $\left( \cup_{i=1}^{l-2} K_b n_w^{i}\right) x_s(1)\subset S_b$. Finally, we conclude that $\partial (S_b) \subset K_b n_w^{-1} \cup K_b n_w^{l-1} \cup K_b h_t(\lambda) \cup K_b h_t(\lambda) ^{-1} \cup K_b x_s(1)$.
\qed
\end{pf}
\subsection{$C_l$}

The Dynkin diagram is the following in this case:

\unitlength 1mm
\begin{picture}(72,12)\label{dynkc}

\put(5,5){ \circle{2} }
\put(5,8){$r_1$ }
\put(7,5){ \line(1,0){8} }

\put(19,5){ \circle{2} }
\put(19,8){$r_2$ }
\put(21,5){ \line(1,0){8} }

\put(33,5){ \circle{2} }
\put(33,8){$r_3$ }
\put(35,5){ \line(1,0){6} }

\put(43,5){ \mbox{\dots }}
\put(56,5){ \circle{2} }
\put(56,8){$r_{l-1}$ }
\put(59,5.5){ \line(1,0){8} }
\put(54,5){ \line(-1,0){6} }

\put(60.1,4.5){\line(1,0){8}}
\put(70,5){ \circle{2} }
\put(70,8){$r_l$ }

\end{picture}
\\It can easily be verified using the Dynkin diagram that $w_{l-1}(r_l)=r_l + 2 r_{l-1}$ and $w_l(r_{l-1} ) =r_{l-1} + r_l$.

Using Lemma \ref{lem} one can see that $w(r_i )= r_{i+1}$ for $i =1,2, \ldots , l-2$.
We also have
\begin{equation*} \begin{split} w(r_{l-1} ) &=w_1 w_2 \ldots w_l (r_{l-1} ) = w_1 w_2 \ldots w_{l-1} (r_l + r_{l-1}) \\ &= w_1 w_2 \ldots w_{l-2} (r_l + r_{l-1}) \mbox{.} \end{split} \end{equation*}
Since $r_l$ is orthogonal to the remaining roots $r_1, r_2, \ldots , r_{l-2}$ we have
\[ w(r_{l-1} ) = r_l + w_1 w_2 \ldots w_{l-2} (r_{l-1} )  \mbox{.}\]
Since $w_i ( r_{i+1} \stp r_{l-1}) =r_i + r_{i+1} \stp r_{l-1}$ for $i=1 \ldots l-2$ we also have
\[ w_1 w_2 \ldots w_{l-2} (r_{l-1} ) = w_1 w_2 \ldots w_{l-3} (r_{l-2} +r_{l-1} ) = r_{1} \stp r_{l-2} +  r_{l-1} \mbox{.}\]
This gives $w(r_{l-1} ) =r_1 +r_2 + \ldots + r_l$.

Using all these observations we can determine a part of the orbit of $\langle w \rangle$ containing $r_1$, which is the following:

\unitlength 1mm
\begin{picture}(101,10)

\put(0,5){\vector(1,0){4}}
\put(5,5){$r_1$ }
\put(7,5){ \vector(1,0){8} }
\put(11,6){\mbox{w}}
\put(17,5){$r_2$ }
\put(19,5){ \vector(1,0){6} }
\put(22,6){\mbox{w}}
\put(26,5){ \mbox{\dots }}
\put(31,5){ \vector(1,0){6} }
\put(34,6){\mbox{w}}
\put(39,5){$r_{l-1}$ }
\put(45,5){ \vector(1,0){8} }
\put(47,6){\mbox{w}}
\put(55,5){$r_l + r_{l-1} + r_{l-2}+ \ldots + r_1 $ }
\put(96,5){ \vector(1,0){4} }

\end{picture}

\begin{lem}
$x_1(1)$, $n_w$ and $h_{r_1} (\lambda )$ generate the Chevalley
group $G$.

\end{lem}
\begin{pf}
The proof can be found in \cite{1}. \qed
\end{pf}

The construction is almost the same as in the case $A_l$.
Let \[ \Gamma _c  = Cay \left( G, \left\{ x_1(1), x_1(-1), n_w,
    n_w^{-1}, h_{r_1}
  (\lambda), h_{r_1} (\lambda) ^{-1} \right\} \right) \mbox{.}\]
 Let
\[ K_c = \langle X_{r_2}, X_{-r_{2} }, X_{r_3},
   X_{-r_{3}}, \ldots ,X_{r_l},  X_{-r_{l} }  \rangle \] and let
 \[ S_c = \cup_{i=0}^{l-2} K_c n_w^{i} \mbox{.} \]

\begin{lem}
  $ \frac{\left| \partial(S_c) \right|}{\left| S_c \right| } \le \frac{6}{l-1}$

\end{lem}
\begin{pf}
Similarly to the previous cases $n_w^{-i} K_c n_w^{i}$ contains $n_w^{-i} X_{r_{i+1}}  n_w^{i} = X_{r_1}$ for $1 \le i \le l-2$ which gives that $n_w^{i}$ is not in $K_c$ if $1 \le i \le l-2$. This proves that $\left| S_c \right| = \left( l-1 \right) \left| K_c \right|$.

Again, $K_c n_w^{i} n_w \subset S_c$ if $1 \le i \le l-3$ and $K_c n_w^{i} n_w^{-1} \subset S_c$ if $i \ne 0$.

It is also easy
to verify that $n_w^{i} x_1( 1)^{\pm 1} n_w^{-i} = \left( x_{i+1} (t) \right)^{ \pm 1 }$ for some $t \in GF(q)^{*}$. Therefore $n_w^{i} x_1( 1)^{\pm 1} n_w^{-i} \in X_{r_{i+1} }$ and $n_w^{i}
h_{r_1}(\lambda )^{\pm 1} n_w^{-i} $ are in the subgroup generated by $X_{r_{i+1} }$ and $X_{-r_{i+1 }}$ for $i=1 \stb l-2$.  Thus the elements
of the right cosets $K_c n_w^{i} x_1 (1)^{\pm 1}$ and $K_c n_w^{i} h_{r_1} (\lambda) ^{\pm 1}$
are in $S_c$ if $1 \le i \le l-2$. This proves that $\partial S_c \subseteq K_c n_w ^{l-1} \cup K_c n_w^{-1} \cup K_c x_1(1)\cup K_c x_1(1)^{-1} \cup K_c h_{r_1}(\lambda) \cup K_c h_{r_1}(\lambda) ^{-1}$, which is the union of $6$ right cosets of $K_c$. Thus $|\partial S_c| \le 6 |K_c|$.
\qed
\end{pf}

\subsection{$D_l$}\label{dn}

The Dynkin diagram in this case is the following:

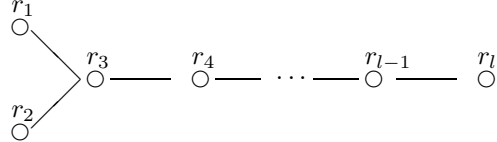
\unitlength 1mm
\begin{picture}(100,21)
\put(9,17){ \circle{2} }
\put(9,19){$r_1$ }
\put(10.5,16.5){ \line(1,-1){6.5} }

\put(9,3){ \circle{2} }
\put(9,5){$r_2$ }
\put(10.5,3.5){ \line(1,1){6.5} }

\put(19,10){ \circle{2} }
\put(19,12){$r_3$ }
\put(21,10){ \line(1,0){8} }

\put(33,10){ \circle{2} }
\put(33,12){$r_4$ }
\put(35,10){ \line(1,0){6} }

\put(43,10){ \mbox{\dots }}
\put(56,10){ \circle{2} }
\put(56,12){$r_{l-1}$ }
\put(59,10){ \line(1,0){8} }

\put(54,10){ \line(-1,0){6} }
\put(71,10){ \circle{2} }
\put(71,12){$r_l$ }

\end{picture}

\begin{lem}
\begin{enumerate}
\item
$x_{r_1}(1)$, $n_w$ and $h_{r_1} (\lambda )$ generate the Chevalley
group $G$ if the rank of $G$ is odd.
\item
$x_{-r_1}(1)$, $x_{r_1}(1)$, $x_3(1)$, $n_w$ and $h_{r_1} (\lambda )$ generate the Chevalley
group $G$ if the rank of $G$ is even.
\end{enumerate}
\end{lem}
\begin{pf}
The proof can be found in \cite{1}.
\qed
\end{pf}
First, we describe a part of the orbit of $\langle w \rangle$ which contains $r_1$.
The root $r_1$ is orthogonal to $r_4, \ldots ,r_l$ hence $w (r_1) =w_1 w_2 w_3 (r_1)$ so we have
\[ w_1 w_2 w_3 (r_1)= w_1 w_2 (r_1 +r_3) = w_1 (r_1 + r_2 +r_3) = -r_1 + r_2 + r_3 +r_1 = r_2 +r_3 \]
and similarly
\[ w (r_2) =w_1 w_2 w_3 (r_2)= w_1 w_2 (r_2 +r_3) = w_1 (-r_2 +r_3 + r_2)  = r_3 + r_1  \mbox{.} \]

Using Lemma \ref{lem}  we get $w(r_i) =r_{i+1}$ for $i=3, \ldots ,l-1$.
This gives that both $w^{i} (r_1)$ and $w^{i} (r_2)$ are of the form
\begin{equation}\label{pot3} r_{i+2} + r_{i+1} + \ldots + r_3 + y \mbox{,}
\end{equation} where $y =r_1$ or $y=r_2$.

\subsubsection{Odd case}
Let us assume that $l$ is odd.

Let \[ \Gamma _d = Cay \left( G, \left\{ x_{r_1}(1), x_{r_1}(1)^{-1}, n_w,
    n_w^{-1}, h_{r_1}
  (\lambda), h_{r_1} (\lambda) ^{-1} \right\} \right) \mbox{.}\]

Let
\[ K_d = \langle X_{r_1}, X_{-r_{1} }, X_{r_2},
   X_{-r_{2} }, \ldots ,X_{r_{l-1} },  X_{-r_{l-1} }  \rangle \] and let
 \[ S_d = \cup_{i=0}^{l-3} K_d n_w^{i} \mbox{.}  \]
\begin{lem}\label{odd2}
  $ \frac{\left| \partial(S_d) \right|}{\left| S_d \right| } \le \frac{2}{l-2}$

\end{lem}
\begin{pf}
It is easy to see that if $0 \le i \le l-3$, then $n_w^{i} x_{r_1} (1)^{\pm 1} n_w^{-i}  =x_{w^{i} (r_1) } (t)^{\pm 1} \in K_d$ for some $t \in GF(q)^*$ since by (\ref{pot3}) the root $w^{i} (r_1)$ is a linear combination with integer coefficients of the fundamental roots $r_1, r_2, \ldots ,r_{l-1}$ and similarly $n_w^{i} h_{r_1} (\lambda)^{\pm 1} n_w^{-i} \in K_d$. It follows that $\partial S_d \subseteq K_d n_w^{l-2} \cup K_d n_w^{-1}$.

It remains to show that for $S_d$ is the union of $l-3$ pairwise disjoint cosets of $K_d$.
Again, $n_w^{i} K_d n_w^{-i}$ contains the subgroup $n_w^{i} X_{r_{l-i} } n_w^{-i} = X_{r_l}$  if $1 \le i \le l-3$ which shows that $n_w^{i} \notin K_d$.
\qed
\end{pf}

\subsubsection{Even case}
Let us assume that $l$ is even.

Let \[ \Gamma _d ' = Cay \left( G, \left\{ x_{-r_1} (\pm 1), x_{-r_1} (\pm 1), x_{r_3}(\pm 1), n_w,
    n_w^{-1}, h_{r_1}
  (\lambda), h_{r_1} (\lambda) ^{-1} \right\} \right) \mbox{.}\]
 Let
\[ K_d ' = \langle X_{r_1}, X_{-r_{1}}, X_{r_2},
   X_{-r_{2} }, \ldots ,X_{r_{l-1}},  X_{-r_{l-1} }  \rangle \] and let
 \[ S_d '= \cup_{i=0}^{l-4} K_d ' n_w^{i} \mbox{.}  \]

\begin{lem}
  $ \frac{\left| \partial(S_d') \right|}{\left| S_d' \right| } \le \frac{2}{l-3}$
\end{lem}
\begin{pf}
It is clear that $w^{i}(-r_1) =-w^{i}(r_1)$ and hence $w^{i}(-r_1)$ is in the root system generated by the roots $r_1, r_2  \stb r_{l-1}$ if $1 \le i  \le l-4$. This shows that $n_w^{i} x_{ \pm r_1} (\pm 1) n_w^{-i} =x_{w^{i} (\pm r_1) } (t) \in  K_d '$ for some $t \in GF(q)^{*}$.

It was proved in Lemma \ref{odd2} that $n_w^{i} h_{r_1} (\lambda)^{\pm 1} n_w^{-i} \in K_d '$ if $0 \le i \le l-4$. Finally, by Lemma \ref{lem} $n_w^{i} x_{r_3} (\pm 1) n_w^{-i} =x_{r_{3+i}}(t)$  for some $t \in K^{*}$ which is in $K_d '$ if $0 \le i  \le l-4$.
It follows that $\partial S_d ' \subseteq K_d ' n_w^{l-3} \cup K_d ' n_w^{-1}$.

It remains to show that $S_d '$ is the union of $l-3$ pairwise disjoint cosets of $K_d '$. This is clear since if $1 \le i \le l-4$, then  $n_w^{i} K_d ' n_w^{-i}$ contains the subgroup $X_{r_l}$ which shows that $n_w^{i} \notin K_d '$.
\qed
\end{pf}

\section{Twisted groups}\label{twisted}

The twisted groups can be obtained as subgroups of Chevalley groups. In order to define twisted groups we need to find a non-trivial symmetry $\rho$  of the Dynkin diagram. We restrict our attention to those twisted groups which are defined using a symmetry of order $2$ and we also assume that the roots in $\Phi$ have the same length. It is well known that such an symmetry $\rho$ can be extended to a unique isometry $\tau$ of $V$ which is the vector space spanned by $\Phi$. We assume that $Aut(GF(q)^{*})$ contains an element of order $2$. Then the Chevalley group $G$ has an automorphism of order $2$, which we denote by $\alpha$ such that $x_r (t)^{\alpha } = x_{\overline{r} } (\overline{k} )$
for every $r \in \pm \Pi$ and $k \in K$, where $\overline{k}  =\tau(k)$ and $\overline{r} = \rho (r)$.

The subgroup $U^{1}$ is the set of elements $u \in U$ such that $u^{\alpha} =u$ and similarly $V^{1}=\left\{ v \in V \mid v^{\alpha } =v \right\}$.
The twisted group $G^{1}$ is generated by $U^{1}$ and $V^{1}$. The subgroups $H^{1}$ and  $N^{1}$ are defined as the intersection of $G^{1}$ with $H$ and $N$, respectively.
We denote by $W^{1}$ the elements $w$ of the Weyl group $W$ such that $ \tau w \tau ^{-1} =w$. There is a natural isomorphism of the group $W^{1}$ to $N^{1}/H^{1}$ and we denote by $n_w^{1}$ the element of $N^{1}  \le N$ which corresponds to $w^{1} \in W^{1}$.

The set of positive roots $\Phi^{+}$ has a partition where the elements of the partition are of the following form:
\begin{flalign*}
 Z&=\left\{ r \mid r \in \Phi^{+} \right\}  \mbox{if } \overline{r} =r \\
Z&=\left\{ r , \overline{r} \mid r \in \Phi^{+} \mbox{and } r+ \overline{r} \notin \Phi \right\} \\
Z&=\left\{ r , \overline{r}, r + \overline{r} \mid r \in \Phi^{+} \mbox{and } r+ \overline{r} \in \Phi \right\} \mbox{.}
\end{flalign*}
We denote by $\Pi^{1}$ the collection of sets which are elements of the partition.
For each set $Z$ in the partition there is a unique element $w_Z \in W^{1}$ which is generated by $\left\{ w_r \mid r \in Z\right\}$ such that $w(Z)= -Z$. These elements are the following:
\begin{flalign*}
 w_Z&=w_r \mbox{ } \mbox{   if  } \mbox{  } Z=\left\{ r \mid r \in \Phi^{+} \right\}  \\
w_Z&=w_r w_{\overline{r} } \mbox{ } \mbox{  if  } \mbox{  }  Z=\left\{ r , \overline{r} \mid r \in \Phi^{+} \mbox{and } r+ \overline{r} \notin \Phi \right\}\\
w_Z&=w_{r + \overline{r} }=w_r w_{\overline{r} } w_r  \mbox{ } \mbox{  if  } \mbox{  } Z=\left\{ r , \overline{r}, r + \overline{r} \mid r \in \Phi^{+} \mbox{and } r+ \overline{r} \in \Phi \right\} \mbox{.}
\end{flalign*}
Every element of $\Pi^{1}$ can be obtained as $w(Z)$, where $w \in W^{1}$ and $Z$ contains a fundamental root. Those sets which contain a fundamental root are called fundamental sets. Moreover, $W^{1}$ is generated by $\left\{ w_Z  \mid Z \in \Pi^{1} \right\}$.

For every $Z \in \Pi^{1}$ we denote by $X_Z$ the subgroup generated by the root subgroups $X_r$ for $r \in Z$ and $X_{Z}^{1} =X_Z  \cap G^1$.

\subsection{$A_{2n-1}^{1}$ }
The fundamental sets in this case are the following:
\[ Z_n =\left\{ r_n \right\}, \mbox{   } Z_i =\left\{ r_i, r_{2n-i} \right\}
\mbox{for  } 1 \le i \le n-1 \mbox{,} \]
and the corresponding elements of the Weyl group $W^{1}$ are:

\[ w_{Z_n} =w_n , \mbox{ and  } w_{Z_i} =w_i w_{2n-i} \mbox{  for  } 1 \le i \le n-1 \mbox{.} \]
We may assume (see \cite[p.233.]{9}) that the subgroups defined above are of the following form :
\begin{flalign*}
X_Z^{1} &= \left\{ x_r(t) \mid t = \overline{t} \right\}  \mbox{ if }Z=\{ r\} \\
   X_Z^{1} &= \left\{ x_r(t) x_{\overline{r}}(\overline{t} )\mid t \in K \right\}  \mbox{ if } Z=\{ r, \overline{r} \} \mbox{.}
 \end{flalign*}

Let $n_w^{1} =n_{w_1^{1}} n_{w_2^{1}} \ldots$ and $h_e= h_{r_1} (t) h_{\overline{r_1} } (\overline{t} )$, where $t$ generates $K^{*}$. In the following in order to simplify notation we write $n_w$ instead of $n_w^{1}$.
\\We also define $x_e =x_{r_1}(1) x_{r_{2n-1} } (1)$ which is an element of $X_{Z_1}^{1}$ and which can also be written
as $x_{r_1}(1) x_{r_1}(1)^{\alpha } = x_{r_1}(1) x_{\overline{r_1 } }(1)$.

\begin{lem}\label{gene}
$x_e$, $n_w$ and $h_e$ generate the group $G^{1}$.
\end{lem}
\begin{pf}
The proof can be found in \cite{1}. \qed
\end{pf}
Let \[ \Gamma _e  = Cay \left( G, \left\{ x_e, x_e ^{-1}, n_w,
    n_w^{-1}, h_e , h_e ^{-1} \right\} \right) \mbox{.}\]
Let
\[ K_e = \langle X_{Z_2}^{1}, X_{-Z_{2} }^{1}, X_{Z_3}^{1},
   X_{-Z_3}^{1}, \ldots ,X_{Z_n }^{1},  X_{-Z_{n} }^{1}  \rangle \] and let
 \[ S_e = \cup_{i=0}^{n-2} K_e n_w^{i} \mbox{.} \]

$K_e$ can be considered as a twisted group which is a subgroup of the Chevalley group generated by the root subgroups $X_{r_2}, X_{-r_2},  \ldots , X_{r_{2n-2}}, X_{-r_{2n-2} }$. The corresponding set of fundamental roots  is $\rho$-invariant and we denote by $\Phi_{2n-3}$ the root system generated by these roots. The restriction of $\rho$ to the set $\left\{ r_2, r_3, \ldots ,r_{2n-2} \right\}$ gives a symmetry of the Dynkin diagram of these roots which extends to an isometry. This isometry is the restriction of $\tau$.
This gives that for $Z \in \Pi^{1}$ the subgroup $X_Z^{1}$ is a subgroup of $K_e$ if and only if $Z \subset \Phi _{2n-3}$.
Clearly, $h_r(t)$ is in $\langle X_Z^{1},X_{-Z}^{1} \rangle \subset G^{1}$ if $Z=\left\{r \right\}$ with $r = \overline{r}$ and if $Z= \left\{ r, \overline{r} \right\}$, then there is homomorphism of $SL_2(K)$ onto $\langle X_Z^{1},X_{-Z}^{1} \rangle \subset G^{1}$ which shows that $x_r(t) x_{\overline{r} } (\overline{t} ) \in G^{1}$ and $h_r(t) h_{\overline{r} } (\overline{t} ) \in G^{1}$.

Conjugating by $n_w^{i} \in N^{1}$ we get the following:
 \begin{equation}\label{konj} \begin{split}n_w^{-i} X_Z^{1} n_w^{i} &= n_w^{-i} (X_Z \cap G^{1}) n_w^{i} = n_w^{-i} X_Z n_w^{i} \cap n_w^{-i} G^{1} n_w^{i} = X_{w^{-i}(Z)} \cap G^{1} \\ &= X_{w^{-i}(Z)}^{1} \mbox{.} \end{split} \end{equation}
\begin{lem}\label{lem15}
  $ \frac{\left| \partial(S_e) \right|}{\left| S_e \right| } \le \frac{6}{n-1}$
\end{lem}
\begin{pf}
We claim that $S_e$ is the union of $n-1$ disjoint subsets. $K_e n_w^{j} =K_e n_w^{j'}$ if and only if $n_w^{j-j'} \in K_e$ so we have to show that $n_w^{i} \notin K_e$ if $1 \le i \le n-2$.
\\We claim that $w^{i} (r_1) =r_{i+1}$ if $1 \le i \le n-2 $.
If $k \le n-3$, then
\[w(r_k)=  w_1 w_{2n-1} \ldots w_{k} w_{2n-k} w_{k+1}   (r_k)
\]
since $r_k$ is orthogonal to $r_j$ if $j \ge k+2$. Therefore
\begin{equation*}\label{pot4} \begin{split} w (r_k) &= w_1 w_{2n-1} \ldots w_{k} (r_k +r_{k+1} ) \\
&= w_1 w_{2n-1} \ldots w_{k-1} (r_{k+1} ) =r_{k+1}  \end{split}
\end{equation*}
 since $r_{k+1}$ is orthogonal to the roots $r_{2n-k} \stb r_{2n-1}$ and  $r_{k+1}$ is orthogonal to $r_1 \stb r_{k-1}$.
It follows that $w^{-i} (Z_{i+1} ) = Z_1$ and hence by equation (\ref{konj}) $X_{Z_{1}}^{1} \subset {n_w}^{-i} S_e n_w^{i} $ if $1 \le i \le n-2$. This proves that $n_w^{i} \notin K_e$ if $1 \le i \le n-2$ hence $\left| S_e \right| = (n-1) \left| K_e \right| $.

It is easy to see that $S_e$ contains $K_e n_w^{i} n_w$ if $i =0, 1, \ldots , n-3$ and $S_e$ contains
$K_e n_w^{i} n_w^{-1}$  if $i = 1, 2, \ldots , n-2$.

We use again the fact that $K_e n_w^{i} g =K_e n_w^{i}$ if and only if $n_w^{i } g n_w^{-i } \in K_e$.
Since $n_w^{i} x_{r_1} (1) n_w^{-i} =x_{w(r_1 ) } (\lambda)$ for some $\lambda \in K$ and $x_e=x_{r_1} (1)x_{r_1} (1)^{\alpha}$ we have
\begin{equation*} \begin{split} &n_w^{i} x_{r_1} (1) x_{r_{1} } (1)^{\alpha} n_w^{-i} = n_w^{i} x_{r_1} (1) n_w^{-i} n_w^{i}x_{r_{1} } (1)^{\alpha} n_w^{-i}\\ = &n_w^{i} x_{r_1} (1) n_w^{-i} (n_w^{i}x_{r_{1} } (1) n_w^{-i})^{\alpha} =x_{w^{i}(r_1 ) } x_{w^{i}(r_1 ) }(\lambda)^{\alpha} = x_{r_{i+1}  } x_{r_{i+1 }   }(\lambda)^{\alpha} \mbox{.} \end{split} \end{equation*}
This shows that $n_w^{i } x_e n_w^{-i } \in X_{Z_{i+1}}^{1}$ which proves that if $i = 1,2, \ldots , n-2$, then $n_w^{i } x_e^{\pm 1} n_w^{-i } \in K_e$ and hence $K_e n_w^{i} x_e^{\pm 1} =K_e n_w^{i}$  since $Z_{i+1} \subset \Phi_{2n-3}$.

We also have $n_w^{i } h_{r_1} (t) h_{\overline{r_1} } (\overline{t} )  n_w^{-i } =h_{r_{i+1} } (\theta) h_{w^{i}  (\overline{r_1} ) }
(\theta' )$ for some $\theta \mbox{,} \theta' \in K$. Using the fact that $w \in W^{1}$ we have $w^{i}  (\overline{r_1} )= \overline{ w^{i}(r_1) }$ so $h_{r_{i+1} } (\theta) h_{w^{i}  (\overline{r_1} ) } (\theta' )= h_{r_{i+1} } (\theta) h_{\overline{
r_{i+1} } } (\theta' )$. Clearly, $n_w^{i } h_e n_w^{-i } \in H^{1}$.
Thus $\theta' =\overline{\theta}$ and $n_w^{i } h_e^{\pm 1}  n_w^{-i } = \left( h_{r_{i+1} } (\theta) h_{\overline{ (r_{i+1} ) } }
(\overline{ \theta ) } \right)^{\pm 1}  \in K_e$ since $r_{i+1} \in \Phi_{2n-3}$ if $i=1, \ldots , n-2$.
This proves that $K_e n_w^{i} h_e^{\pm 1} =K_e n_w^{i}$ if $i=1, \ldots , n-2$ and hence $\partial S_e \subset K_e n_w^{n-1} \cup K_e n_w^{-1} \cup K_e x_e \cup K_e x_e^{-1} \cup K_e h_e \cup K_e h_e^{-1} $, finishing the proof of Lemma \ref{lem15}.
\qed
\end{pf}

\subsection{$D_n^{1}$}
The fundamental sets in this case are the following:
\[ Z_1 =\left\{ r_1, r_2 \right\}, \mbox{   } Z_i =\left\{ r_{i+1}  \right\}
\mbox{for  } 2 \le i \le n-1 \mbox{,}
\]
and the corresponding elements of the Weyl group $W^{1}$ are:
\[ w_{Z_1} =w_1 w_2 , \mbox{ and  } w_{Z_i} =w_{i+1} \mbox{  for  } 2 \le i \le n-1 \mbox{.} \]
Let $n_w =n_{w_1^{1}} n_{w_2^{1}} \ldots n_{w_{n-1}^{1}}$ and $h_f= h_{r_1} (t) h_{\overline{r_1} } (\overline{t} )$, where $t$ generates $K^{*}$.
\\We also define $x_f =x_{r_1}(1) x_{r_{2} } (1)$ which can also written as $x_{r_1}(1) x_{r_1}(1)^{\alpha } = x_{r_1}(1)
x_{\overline{r_1 } }  (1  )$.

\begin{lem}\label{genf}
$x_f$, $n_w$ and $h_f$ generate the group $G^{1}$.
\end{lem}
\begin{pf}
The proof can be found in \cite{1}. \qed
\end{pf}
Let \[ \Gamma _f   = Cay \left( G, \left\{ x_f, x_f ^{-1}, n_w,
    n_w^{-1}, h_f , h_f ^{-1} \right\} \right) \mbox{.}\]
Let
\[ K_f = \langle X_{Z_1}^{1}, X_{-Z_{1} }^{1}, X_{Z_2}^{1},
   X_{-Z_2}^{1}, \ldots ,X_{Z_{n-2} }^{1},  X_{ -Z_{n-2} }^{1}  \rangle \] and let
 \[ S_f = \cup_{i=0}^{n-3} K_f n_w^{i} \mbox{.} \]

We denote by $\Phi_{n-1}$ the root system generated by the fundamental roots $r_1, r_2, \ldots ,r_{n-1}$.
 \begin{lem}
  $ \frac{\left| \partial(S_f) \right| }{\left| S_f \right| } \le \frac{2}{n-2}$
\end{lem}
\begin{pf}
The Coxeter element in this case is exactly the same as in subsection \ref{dn}. This gives that $n_w^{i} (r_{n-i}) =r_n$ for $0 \le i \le n-3$.
The fundamental sets $Z_2, Z_3, \ldots ,Z_{n-1}$ consist of only one element thus $n_w^{i} S_f n_w^{-i}$ contains
 $X^{1} _{w^{i} ( Z_{n-1-i} ) } = X_{w^{i} ({r_{n-i} ) }  } = X_{r_n} = X_{Z_{n-1} }^{1} $ if $1 \le i \le n-3$ since $S_f$ contains $X_{n-i}^{1}$. This proves that if $1 \le i \le n-3$, then $n_w ^{i} \notin K_f$. Thus $S_f$ is the union of $n-2$ disjoint subsets of the same cardinality. Therefore $\left| S_f \right|= (n-2) \left| K_f \right| $.

Using the definiton of $S_f$ one can see that $K_f n_{w}^{i}  n_w \subset S_f$ if $i= 0,1, \ldots ,n-4$ and
$K_f n_w^{i} n_w^{-1} \subset S_f$ if $i= 1, \ldots ,n-3$.

The elements $n_w^{i} x_f n_w ^{-i}$ are of the form $x_r (t) x_{\overline{r} } (\pm \overline{t} )$ for some $r \in \Phi$ and $t \in K^{*}$.
In order to prove that these elements are in $K_f$ for $i=0, 1, \ldots ,n-3$ we only have to show that $r \in \Phi_{n-1}$. Using the fact that the Coxeter element in this case is the same as in Section \ref{dn} we have that both $w^{i} (r_1)$ and $w^{i} (r_2)$ are of the form $r_1 + r_3 + r_4 +\ldots + r_{i+1}$ or $r_2 + r_3 + r_4 +\ldots + r_{i+1}$.
These roots are clearly in the root system generated by the fundamental roots $r_1, r_2, \ldots ,r_{l-1}$ if $i \le n-2$. This proves that $n_w^{i} x_f^{\pm} n_w ^{-i}$ is in $K_f$ if $0 \le i \le n-3$ and hence $S_f x_f^{\pm } \subset S_f$.

Similarly, the elements $n_w^{i} h_f n_w ^{-i}$ are of the form $h_r (t) h_{\overline{r} } ( \overline{t} )$ for some $r \in \Phi$ and $t \in K^{*}$ and it is easy to see that $r \in \Phi_{n-1}$ if $0 \le i \le n-3$. This proves that $n_w^{i} h_f^{\pm} n_w ^{-i}$ is in $K_f$ if $0 \le i \le n-3$ and hence $S_f h_f^{\pm } \subset S_f$.
\qed
\end{pf}

 \subsection{$A_{2n}^{1}$ }
The fundamental sets are the following:
\[ Z_1 =\left\{ r_n, r_{n+1}, r_n + r_{n+1} \right\}, \mbox{   } Z_i =\left\{ r_{n+1-i}, r_{n+i} \right\}
\mbox{ for  } 2 \le i \le n \mbox{.} \]

Let $n_w^{1} =n_{w_1^{1}} n_{w_2^{1}} \ldots n_{w_n ^{1}}$ and $h_g= h_{r_n} (t) h_{\overline{r_n} } (\overline{t} )$, where $t$ generates $K^{*}$.
\\We also define $x_g =x_{r_n}(1) x_{r_{n+1} }(1) x_{r_n +r_{n+1} }  (k)$  with $k + \overline{k} =1$.

\begin{lem}\label{geng}
$x_g$, $n_w$ and $h_g$ generate the group $G^{1}$.
\end{lem}
\begin{pf}
The proof can be found in \cite{1}.
\qed
\end{pf}

Let \[ \Gamma _g  = Cay \left( G, \left\{ x_g, x_g ^{-1}, n_w,
    n_w^{-1}, h_g , h_g ^{-1} \right\} \right) \mbox{.}\]
Let
\[ K_g = \langle X_{Z_1}^{1}, X_{-Z_{1} }^{1}, X_{Z_2}^{1},
   X_{-Z_2}^{1}, \ldots ,X_{Z_{n-1} }^{1},  X_{-Z_{n-1} }^{1}  \rangle \] and let
 \[ S_g = \cup_{i=0}^{n-2} K_g n_w^{i} \mbox{.} \]

\begin{lem}
  $ \frac{\left| \partial(S_g) \right| }{\left| S_g \right| } \le \frac{2}{n-1}$
\end{lem}
\begin{pf}
First, we show that $S_g$ is the union of $n-1$ disjoint subsets of the same cardinality.
It is enough to show that $n_w^{i} \notin K_g$ for $i=1, \ldots , n-2$. This will be done by proving that $X_{Z_{n} }^{1}$
is contained in $n_w^{i} K_g n_w^{-i}$. Using equation (\ref{konj}) we only have to show that $w^{i} (Z_{n-i} )= Z_{n}$ for $i=1, \ldots n-2$.

The fundamental root $r_{k+1}$ is contained in $Z_{n-k}$. Let us assume that $1 \le k \le n-2$.
\begin{flalign*} w(r_{k+1} )&= w_n w_{n+1} w_n w_{n-1} w_{n+2} \ldots w_1 w_{2n} (r_{k+1} ) \\ &= w_n w_{n+1} w_n w_{n-1} w_{n+2} \ldots w_{k+1} w_{2n-k} w_k (r_{k+1} ) \end{flalign*}
since $r_{k+1}$ is orthogonal to the roots $r_j$ if $j>n$ or $j<k-1$.
Clearly, $w_{k+1} w_{2n-k} w_k (r_{k+1} ) =r_k$ so
\[  w(r_{k+1} ) =w_n w_{n+1} w_n \ldots w_{k+2} w_{2n-k-1} (r_k) =r_k \]
since the remaining reflections fix $r_k$.

One can see by induction that $w^{i}(r_{i+1}) =r_1$ for $i=1, \ldots n-2$ and since $w \in W^{1}$  we have $w^{i} (\overline{r_{i+1} }) =\overline {w^{i} (r_{i+1} ) }=r_{2n}$ and hence $w^{i} (Z_{n-i} ) =Z_n$. This proves that for $i=1, \ldots  ,n-2$ the subgroup $n_w^{i} (K_g) n_w^{-i}$ contains  $X_{Z_n}^{1}$. Therefore $\left| S_g \right| = (n-1) \left| K_g \right|$.

The definition of $S_g$ shows that $K_g n_w^{i} n_w \subset S_g$ if $i \ne n-2$ and $K_g n_w^{i} n_w^{-1} \subset S_g$ if $i \ne 0$.
It remains to investigate the elements of the form $n_w^{i} x_g^{\pm} n_w^{-i}$ and $n_w^{i} h_g^{\pm} n_w^{-i}$.

We claim, that $w^{i} (r_n) =r_n+ r_{n-1} + \ldots +r_{n-i}$ if $i \le n-2$.
Using the orthogonality of the fundamental vectors $r_j$, $r_{k}$, where $\left| j-k \right| \ge 2$ we get the following:
\begin{equation}\label{eg4} \begin{split}   w(r_n) &=w_n w_{n+1} w_n w_{n-1} w_{n+2}\ldots w_1 w_{2n}(r_n) \\&= w_n w_{n+1} w_n w_{n-1} (r_n) =w_n w_{n+1} (r_{n-1})=r_{n-1} +r_n \mbox{.} \end{split} \end{equation}
Similarly, if $1 \le k \le n-2$, then
\begin{equation}\label{eg3} \begin{split} w(r_{n-k} ) &=w_n w_{n+1} w_n w_{n-1} w_{n+2}\ldots w_1 w_{2n}(r_{n-k} ) \\&= w_n w_{n+1} w_n \ldots w_{n-k} w_{n+k+1} w_{n-k-1} (r_{n-k}) \\&=w_n w_{n+1} w_n \dots  w_{n-k+1} (r_{n-k-1})=r_{n-k-1} \mbox{.}\end{split} \end{equation}
Since $w$ is linear we get using (\ref{eg4}) and (\ref{eg3}) that \begin{equation}\label{eg5} w^{i} (r_n) =r_n+ r_{n-1} + \ldots +r_{n-i} \mbox{.} \end{equation}
By observing equations (\ref{eg5}) one can see that if $i=0, \ldots ,n-2$, then both $r_1$ and $r_{2n}$ are orthogonal to $w^{i} (r_n)$ and similarly $r_1$ and $r_{2n}$ are orthogonal to $w^{i}(r_{n+1} ) = w^{i}(\overline{r_n}) = \overline{w^{i}(r_n)} =r_{n+i+1}$.
 This shows that for $w'=w_n w_{n+1} w_n w_{n-1} w_{n+2}\ldots w_2 w_{2n-1}$ we have $w^{i}(r_n) =(w')^{i}(r_n)$ and $w^{i}(r_{n+1} ) =(w')^{i}(r_{n+1} )$. Therefore $w^{i}(r_n + r_{n+1} ) =(w')^{i}(r_n + r_{n+1} )$.
Moreover, $n_w^{i} x_g n_w^{-i} =n_{w'}^{i} x_g n_{w'}^{-i}$ and $n_w^{i} h_g n_w^{-i} =n_{w'}^{i} h_g n_{w'}^{-i}$.

Clearly, $n_{w'} \in K_g$ and hence the elements $n_{w'}^{i} x_g^{\pm} n_{w'}^{-i}$ and $n_{w'}^{i} h_g^{\pm} n_{w'} ^{-i}$ are in $K_g$ if $i =0, 1, \ldots n-2$.
\qed
\end{pf}

In order to finish the proof of Theorem (\ref{thm1}) we have to verify that the for those sets $S$ for which boundary $\partial(S)$ is relatively small we have $\left| S \right| \le \frac{\left| G \right|}{2}$.
The order of the investigated simple groups is the following:
\[ \begin{matrix}
A_l(q): & \frac{1}{(n+1,q-1)} q^{\frac{n(n-1)}{2} } \prod_{i=1}^{n} \left( q^{i+1} -1 \right) \\
B_l(q): & \frac{1}{(2,q-1)} q^{n^{2} } \prod_{i=1}^{n} \left( q^{2i} -1 \right) \\
C_l(q): & \frac{1}{(2,q-1)} q^{n^{2} } \prod_{i=1}^{n} \left( q^{2i} -1 \right) \\
D_l(q): & \frac{1}{(4,q^{n}-1)} q^{\frac{n(n-1)}{2} } \prod_{i=1}^{n} \left( q^{2i} -1 \right) \\
A_l(q^{2}) ^{1}: & \frac{1}{(n+1,q+1)} q^{\frac{n(n-1)}{2} } \prod_{i=1}^{n} \left( q^{i+1} - \left( -1 \right)^{i+1} \right) \\
D_l(q^{2}) ^{1}: & \frac{1}{(4,q^{n}+1)} q^{n(n-1)} \left( q^{n} +1 \right) \prod_{i=1}^{n} \left( q^{2i}  -1  \right)
 \end{matrix} \]
It is easy to see that such a simple group can not have a subgroup of index at most $2l$, finishing the proof of Theorem \ref{thm1}.

\section{Identification}\label{identification}
In this section we give explicit generators of the Cayley graphs that we investigated in Section \ref{chevalley} and \ref{twisted}. We also show how to find the subsets of the vertices $S$ for which $\partial S$ is relatively small. We only handle the case of Special Linear Groups which can easily be transformed to the case of the Projective Special Linear Groups which is clearly the easiest one.
This example includes the original idea which was extended to several different series of simple groups. In order to show the simplicity of the original construction we forget about the machinery which was built up before.

Let
\[ A_l =
\begin{pmatrix}
1 & 1 \\
& 1 \\
& & 1  \\

& & & \ddots \\
& & & & 1
\end{pmatrix} \mbox{,}
\]
where $A_l \in GF(q)^{(l+1) \times (l+1)}$.
Let
\[ B_l=
\begin{pmatrix}
0 & 0 & 0 &\dots & (-1)^{l}\\
1 & 0 \\
0& 1 & 0 &  \\
0& 0&  \ddots & \ddots \\
0& &\dots & 1 & 0
\end{pmatrix} \mbox{.}
\]
We denote by $C_l $ the diagonal matrix  $diag(\frac{1}{\lambda},
\lambda, 1, 1, \ldots ,1) \in GF(q)^{(l+1) \times (l+1)}$, where $\lambda$ generates $GF(q)^{*}$.

We denote by $e_{i,j}$ the matrix with $1$ in the $(i,j)$-th
position and zeros everywhere else and let $T_{i,j} (\delta)= I +
\delta e_{i,j} $, where $I$ denotes the identity matrix. Using this notation we can write $A_l=T_{1,2} (1)$.

The standard generator $x_{r_1}(1)$ of the Chevalley group given in Subsection \ref{A_l} corresponds to the matrix $A_l$ and the Coxeter element $n_w$ can be identified with $B_l$. Finally, $C_l$ plays the role of $h_{r_1} (\lambda)$.

Clearly, $T_{i,j}(\alpha) T_{i,j}(\beta) = T_{i,j}(\alpha + \beta)$ and
$\left[ T_{i,j}(\alpha), T_{j,k} (\beta) \right] =T_{i,k}(\alpha \beta)
$ if $i \ne k$, where $\left[ g,h \right] =g^{-1} h^{-1}gh$ denotes the
commutator of $g$ and $h$.

\begin{lem}\label{pr17}
For every $l \in \mathbb{N}$ the set $\left\{ A_l, B_l, C_l \right\}$ forms
a generating set of $SL(l+1,q)$.
\end{lem}
\begin{pf}
We fix the size of the matrices and hence we can write $A= A_l$,
$B=B_l$ and $C=C_l$. Let $H= \langle A, B, C \rangle$. It is enough to
verify that $T_{i,j}(\delta) \in H$ for every $i \ne j$ and $\delta
\in GF(q)$.

It is easy to see that $A^{C^k} =T_{1,2}(1)^{C^k} =T_{1,2} (\lambda
^{2k})$. Using $T_{1,2}(\mu) T_{1,2}(\eta) =T_{1,2}(\mu \eta)$ we get that $T_{1,2}(\delta) \in \langle A, C \rangle \le H$ for
every $\delta \in GF(q)$.
For $i \ne j $ we have $B^{k} T_{i,j}(\delta) B^{-k} =T_{i+k, j+k}(\pm \delta)$, where the indices are taken modulo $l+1$
and hence $T_{i,i+1} (\delta) \in H$ for every $1 \le i \le l$ and for
 every $\delta \in GF(q)$. This implies that for every $1 < l \le l+1$  and for
every $\delta \in GF(q)$
\[ \left[ \ldots \left[
   \left[ T_{1,2}(\delta), T_{2,3}(1) \right], T_{3,4} (1) \right]
\ldots ,T_{k-1,k}(1) \right] =T_{1,k}(\delta) \in H \mbox{.} \]
Using again the fact that $B^{k} T_{1,l}(\delta) B^{-k}  =
T_{1+k, l+k}(\pm \delta)$ we get that $T_{i,j} (\delta) \in H$ for every
$i \ne j$ and for every $\delta \in GF(q)$.
\qed
\end{pf}
Let
\[
S_0 = \left\{ \begin{pmatrix} D & 0 \\ 0 & E  \end{pmatrix} \in
  SL(l+1,q) \Bigg| D \in SL(l,q), E \in
  SL(1,q)  \right\} \mbox{.}
\]
For every $1 \le i \le l$ we define
\[ S_i =  S_0 B^i  \mbox{.}
\]
Finally, let
\[ S = \bigcup _{i=0}^{l-1} S_i  \mbox{.} \]
 It is easy to see that $\left| S
   \right| < \frac{\left| SL(l+1,q) \right| }{2}$ if $l \ge 1$.
\begin{lem}\label{pr18}
$ \frac{\left| \partial(S) \right|}{\left| S \right| } \le \frac{6}{l}$

\end{lem}
\begin{pf}
Every element of $S$ has exactly $l$ columns with
$0$ in the last row, exactly $1$ column with $0$ in the first $l$ and $1$ in the last row. The
sets $ S_i$ are pairwise disjoint since an invertible
matrix can not have a column with only zero entries. Furthermore, they all
have the same cardinality since $S_0$ is a subgroup of $SL(n,q)$
and $S_i$ are right cosets of $S_0$ in $SL(n,q)$.

It is easy to see that $S B \setminus S \subseteq S_0 B^l =S_l$ and $S B^{-1} \setminus S \subseteq S_0 B^{-1}$.
The remaining elements of $\partial S$ are of the form $MA$, $MC$ and $MA^{-1}$, $MC^{-1}$ where $M \in S$.

Let us assume that  $M \in S_i$. Then
\[ M= \begin{pmatrix} D & 0 & D' \\ 0 & 1 & 0  \end{pmatrix}
   \]
for some $D \in GF(q)^{l,l-i}$ and $D' \in GF(q)^{l,i}$.
Multiplying a matrix $M$ by $A$
or $A^{-1}$ from the right only modifies the second column of
$M$. Therefore if $M \in S_i$ with $i \ne l,l-1$, then it is easy to see that $MA, MA^{-1}
\in S_i$.

Multiplying a matrix $M$ by $C$ or $C^{-1}$ from the right only modifies the first and the second columns of $M$ thus if $M \in S_i$ with $i
\ne l,l-1$, then $M C^{\pm 1} \in S_i$.

This gives that $\partial S \subseteq S_l \cup S_0 B^{-1} \cup S_{l-1} A \cup S_{l-1} A^{-1}
\cup S_{l-1} C \cup S_{l-1} C^{-1}$ since $S =\cup _{i=0}^{l-1} S_i$.

\qed
\end{pf}
\section*{Acknowledgement}
The author is grateful to L\'aszl\'o Pyber for many valuable suggestions during the research.

\end{document}